\documentclass[12pt]{article}
\usepackage{amsmath}
\usepackage{theorem}
\usepackage{amssymb}

\theorembodyfont{\sl}
\newtheorem{thm}{Theorem}[section]

\newtheorem{lem}[thm]{Lemma}

\newtheorem{prop}[thm]{Proposition}

\newtheorem{definition}[thm]{Definition}
\newenvironment{de}{\begin{definition}\rm}{\end{definition}}
\newtheorem{example}[thm]{Example}
\newenvironment{exam}{\begin{example}\rm}{\end{example}}
\newtheorem{remark}[thm]{Remark}
\newenvironment{rem}{\begin{remark}\rm}{\end{remark}}

\newcommand{\openbox}{\leavevmode
  \hbox to.77778em{%
    \hfil\vrule
  \vbox to.675em{\hrule width.6em\vfil\hrule}%
  \vrule\hfil}} \newcommand{\proofname}{Proof}
\newenvironment{proof}[1][\proofname]{\par\normalfont
  \trivlist\item[\hskip\labelsep\itshape #1:]\ignorespaces
  }{\hspace*{1cm}\hspace*{\fill}\openbox \medskip\endtrivlist}

\newcommand{\eqr}[1]{~\mbox{$(${\rm \ref{#1}}$)$}}

\newcommand{\Smallmat}[2]{\mbox{%
\footnotesize{$\begin{bmatrix}{#1}\\{#2}\end{bmatrix}$}}}
\newcommand{\Smallpmat}[2]{\mbox{%
\footnotesize{$\begin{pmatrix}{#1}\\{#2}\end{pmatrix}$}}}

\newcommand{\dual}[2]{\mbox{$\langle{#1},{#2}\rangle$}}
\newcommand{\ddual}[2]{\mbox{$\langle\!\langle{#1},{#2}\rangle\!\rangle$}}

\newcommand{\F}{{\mathbb F}}
\newcommand{\N}{{\mathbb N}}
\newcommand{\Z}{{\mathbb Z}}
\newcommand{\R}{{\mathbb R}}
\newcommand{\A}{{\mathcal A}}
\newcommand{\B}{{\mathcal B}}
\newcommand{\C}{{\mathcal C}}
\newcommand{\D}{{\mathcal D}}
\newcommand{\Scal}{{\mathcal S}}
\newcommand{\im}{{\rm im}}
\newcommand{\rank}{{\rm rank}\,}

\newcommand{\vsp}{\vspace*{-1mm}}
\newcommand{\junk}[1]{}

\newcommand{\lf}{\left\lfloor}
\newcommand{\rf}{\right\rfloor}
\def\2stack#1#2{\mathrel{\mathop{#1}\limits_{#2}}}

\title{Duality between Multidimensional Convolutional Codes and 
  Systems\footnote{Supported in part by NSF grant DMS-96-10389.
    Paul Weiner would like to thank the Center for Applied
    Mathematics at Notre Dame for a fellowship which financially
    supported the presented research.}}

\author{
  Heide Gluesing-Luerssen\\
  {\small Department of Mathematics\vsp}\\
  {\small University of Notre Dame\vsp}\\
  {\small Notre Dame, Indiana 46556-5683, USA\vsp}\\
  {\small {\em e-mail:} gluesing-luerssen.1@nd.edu}\\
  {\small {\em or\/} gluesing@mathematik.uni-oldenburg.de} \and
  Joachim Rosenthal\\
  {\small Department of Mathematics\vsp}\\
  {\small University of Notre Dame\vsp}\\
  {\small Notre Dame, Indiana 46556-5683, USA\vsp}\\
  {\small {\em e-mail:} Rosenthal.1@nd.edu} \and
  Paul A.~Weiner\\
  {\small Department of Mathematics\vsp}\\
  {\small Saint Mary's University of Minnesota\vsp}\\
  {\small Winona, MN 55987, USA\vsp}\\
  {\small {\em e-mail:} pweiner@smumn.edu}}

\begin{document}\maketitle

\begin{center}
  {\em Dedicated to Diederich Hinrichsen on the occasion of his
    60th birthday}
\end{center}
\section{Introduction}
\setcounter{equation}{0}

Data transmission over noisy channels requires implementation of
good coding devices. Convolutional codes belong to the most
widely implemented codes. These codes represent in essence
discrete time linear systems over a fixed finite field $\F$.
Because of this reason a study of convolutional codes requires a
good understanding of techniques from linear systems theory.

Multidimensional convolutional codes generalize (one dimensional)
convolutional codes and they correspond to multidimensional
systems widely studied in the systems literature.
(See~\cite{ob90} and its references). These codes are very
suitable if e.g. the data transmission requires the encoding of a
sequence of pictures and we will explain this at the end of this
section.

In the sequel we will assume that a certain message source is
already encoded through a sequence of vectors $m_i\in\F^k$,
$i=1,\ldots,\gamma$. If every vector in $\F^k$ is a valid message
word, then the change of one coordinate of a vector $m\in\F^k$
will result in another valid message vector $\tilde{m}\in\F^k$ and 
the error can neither be detected nor corrected.
In order to
overcome this difficulty one can add some redundancy by constructing 
an injective linear map
$$
\varphi:\F^k\longrightarrow \F^n
$$
having the property that the Hamming distance
dist\,$\left(\varphi(m_1),\varphi(m_2)\right)$, that is the
number of different entries in the vectors  $\varphi(m_1)$ and 
$\varphi(m_2)$, is at least $d$ whenever
$m_1\neq m_2$. If one transmits the $n$-vector $\varphi(m_1)$
instead of the $k$-vector $m_1$ then it is possible to correct up
to $\lf \frac{d-1}{2}\rf$ errors for every transmitted
$n$-vector, for details see Lemma~\ref{L0248}.

The linear transformation $\varphi$ defines an encoder and
$\im(\varphi)\subset \F^n$ is called a linear block code. In order
to describe the encoding of a whole sequence of message words
$m_0,m_1,\ldots,m_\gamma\in\F^k$ it will be convenient to
introduce the polynomial vector $m(z):=\sum_{i=1}^\gamma
m_iz^i\in\F^k[z]$. The encoding procedure is then compactly
written by:
$$
\hat{\varphi}:\F^k[z]\longrightarrow \F^n[z],\ \ m(z)\longmapsto
\hat{\varphi}(m(z))= \sum_{i=1}^\gamma \varphi(m_i)z^i.
$$

If $\D$ denotes the polynomial ring $\D=\F[z]$ then one
immediately verifies that $\hat{\varphi}$ describes an injective
module homomorphism between the free modules $\D^k$ and $\D^n$
and $\im\left( \hat{\varphi}(\D^k)\right)\subset\D^n$ is a
submodule.

In general not every injective module homomorphism between $\D^k$
and $\D^n$ is of this form. Indeed $\hat{\varphi}$ has the
peculiar property that the $i$-th term of $\hat{\varphi}(m(z))$
only depends on the $i$-th term of $m(z)$. In other words the
encoder $\hat{\varphi}$ has `no memory'. In general it is highly
desirable to invoke encoding schemes where
$\hat{\varphi}:\D^k\rightarrow\D^n$ is an arbitrary injective
module homomorphism. The image of such a module homomorphism is
then called a ($1$-D) convolutional code.

$1$-D convolutional codes are very much suited in the encoding of
sequences of message blocks. Sometimes it might be desirable that
the data is represented through polynomial rings in several
variables. This leads us then to the definition of a
$m$-dimensional convolutional code whose basic properties we
intend to study in this paper.

The following example will illustrate the usefulness of
multidimensional convolutional codes.

\newpage
\begin{exam}
  Let $\D=\F[z_1,z_2,z_3]$ be the polynomial ring in the indeterminates
  $z_1,z_2,z_3$. A whole motion picture (without sound) can be described
  by one element of $\D^k$. Indeed if $f\in\D^k$,
  $$
  f(z_1,z_2,z_3)=\sum_{x=0}^\xi\sum_{y=0}^\rho\sum_{t=0}^\tau
  f_{(x,y,t)}z_1^xz_2^yz_3^t\in \F^k[z_1,z_2,z_3]
  $$
  then we can view the vector $f_{(x,y,t)}\in\F^k$ as
  describing the color and the intensity of a pixel point with
  coordinates $(x,y)$ at time $t$.
  
  In practice the encoding of the element $f(z_1,z_2,z_3)\in\D^k$ is
  done in the following way. At a particular time instance $t$
  all the data vectors $f_{(x,y,t)}$ are combined into a large
  vector $\hat{f}_t\in\F^K$, where $K$ depends on the size of $k$
  and the number of pixel points on the screen. In this way we
  can identify each element $f(z_1,z_2,z_3)\in\D^k$ of above type with
  a polynomial vector $\hat{f}(z_3)\in\F^K[z_3]$. The vector
  $\hat{f}(z_3)$ is then encoded with a usual $1$-D encoding
  scheme. This encoding scheme is shift invariant with respect to
  time but it is in general not shift invariant with respect to
  the $z_1$ and $z_2$ directions on the screen.
  
  In order to achieve an encoding scheme which is also
  shift-invariant with respect to the coordinate axes of the
  screen one can do the following: Construct an injective module
  homomorphism $\hat{\varphi}:\D^k\rightarrow\D^n$. The image
  then describes a $3$-dimensional convolutional code which is
  invariant with respect to time and both coordinate axes. The
  transmission of an element $\D^n$ is then done by choosing a
  term order among the monomials of the form $z_1^xz_2^yz_3^t$.
\end{exam}

\section{Multidimensional Convolutional Codes}
\setcounter{equation}{0}

In this section we introduce multidimensional convolutional codes
as submodules of $\D^n$, where $\D$ denotes a polynomial ring in
$m$ variables. Our presentation in this section follows
closely~\cite[Chapter~2]{we98t}.

We begin by setting some notations.  Let $\F$ be any finite field
and define $\D=\F[z_1,\ldots,z_m]$ to be the polynomial ring in
$m$ indeterminates over $\F$.  We will mainly use the shorter
form
\[
\D=\F[z]=\Big\{\sideset{}{'}\sum_{\alpha\in\N^m}
f_{\alpha}z^{\alpha}\,\Big|\, f_{\alpha}\in\F\Big\},
\]
where for $\alpha=(\alpha_1,\ldots,\alpha_m)\in\N^m$ the notation
$z^{\alpha}$ stands for $z_1^{\alpha_1}\cdot\ldots\cdot
z_m^{\alpha_m}$ and where $\sum'$ means this sum being finite.
Note that $\D$ is $\F$-isomorphic to the $m$-dimensional finite
sequence space
\[
\Scal=\{f:\N^m \longrightarrow \F \mid f \text{ has finite
  support}\},
\]
the isomorphism given by
\[
\begin{split}
  \psi\,:\ \Scal & \longrightarrow  \D \\
  f & \longmapsto \sum_{\alpha\in\N^m} f(\alpha) z^{\alpha}
\end{split}
\]
One can visualize the elements of $\Scal$ by using the integer
lattice of the first quadrant of $\R^m$ and attaching the element
$f(\alpha_1,\ldots,\alpha_m)\in\F$ to the point with coordinates
$(\alpha_1,\ldots,\alpha_m)$.  It is convenient to omit the
attachment if $f(\alpha_1,\ldots,\alpha_m)=0$.

\begin{exam}  \label{L0048}    
  We visualize the polynomial
  $f(z_1,z_2)=1+2z_1^2+2z_1z_2\in\F_3[z_1,z_2]$ as well as $z_1
  f(z_1,z_2)$.
  $$
  \begin{array}{ccc}
   \begin{array}{llll}
      \bullet & \bullet _2 & \bullet & \bullet \\
      \bullet _1 & \bullet & \bullet _2 & \bullet
    \end{array} & 
  \stackrel{\psi}{\longrightarrow} & 1+2{z_1}^2+2z_1 z_2 \\[.15in]
   & & \downarrow z_1 \cdot \\[.15in]
  \begin{array}{llll}
    \bullet & \bullet & \bullet _2 & \bullet \\
    \bullet & \bullet _1 & \bullet & \bullet _2
   \end{array} &
  \stackrel{\hskip.05in \psi^{-1}}{\longleftarrow}
     & z_1+2{z_1}^3+2{z_1}^2  z_2 \,\,.
  \end{array}
  $$
\end{exam}

As the example indicates, multiplication with $z_i$ in the ring
$\D$ corresponds to the forward shift along the $i$th axis in
$\Scal$.  This can be verified with the help of the following
commutative diagram\medskip

\unitlength1.1cm \mbox{}\hspace{3cm}%
\begin{picture}(10,4)
  \put(2,3){$\Scal$} \put(2.3,3.1){\vector(1,0){1.6}}
  \put(4,3){$\D$} \put(2.1,2.9){\vector(0,-1){1.6}}
  \put(4.1,2.9){\vector(0,-1){1.6}} \put(2,1){$\Scal$}
  \put(3.9,1.1){\vector(-1,0){1.6}} \put(4,1){$\D$}
  \put(0.3,3.6){\mbox{$\big(f(\alpha)\big)_{\alpha\in\N^m}$}}
  \put(4.5,3.6){$\sum_{\alpha\in\N^m}f(\alpha)z^{\alpha}$}
  \put(-0.2,0.4){\mbox{$\big(f(\alpha-e_i)\big)_{\2stack{\alpha
  \in\N^m}{\!\!\!\!\!\!\alpha_i\not=0}}$}}
  \put(4.5,0.4){\mbox{$\sum_{\2stack{\alpha
  \in\N^m}{\!\!\!\!\!\!\alpha_i\not=0}}f(\alpha-e_i)z^{\alpha}$}}
  \put(0.7,3.3){\vector(0,-1){2.6}}
  \put(0.65,3.3){\line(1,0){.1}}
  \put(5.5,3.3){\vector(0,-1){2.6}}
  \put(5.45,3.3){\line(1,0){.1}} \put(2.1,3.7){\vector(1,0){2.1}}
  \put(2.1,3.65){\line(0,1){.1}}
  \put(4.3,0.5){\vector(-1,0){1.9}}
  \put(4.3,0.45){\line(0,1){.1}} \put(3,3.2){\mbox{\footnotesize
      $\psi$}} \put(3,1.2){\mbox{\footnotesize $\psi^{-1}$}}
  \put(4.2,2){$z_i \cdot$}
\end{picture}\medskip

Here $e_i\in\N^m$ denotes the $i$th standard basis vector.

Throughout this paper a code is defined to be an $\F$-linear
subspace of some $\Scal^n$ which is invariant under the forward
shifts along all axes.  By virtue of the above diagram this can
simply be phrased as

\begin{de}\label{D-code}
  A linear $m$-dimensional  convolutional code (for short, $m$-D code) 
  of length $n$ over $\F$ is a
  $\D$-submodule of $\D^n$.  An element of a code is said to be a
  codeword.
\end{de}

\begin{rem}                      \label{rem1}
  In the coding literature (see e.g.~\cite{pi88}) convolutional
  codes are usually not restricted to sequence spaces whose
  elements have finite support. There is however no engineering
  reason behind this.  After all every transmitted message
  created by mankind did have finite length. Convolutional codes
  with finite support were first studied by Fornasini and
  Valcher~\cite{fo94a1,fo98a,va94}. These authors did define a
  convolutional code as a submodule of $\tilde{\D}^n$ where
  $\tilde{\D}$ represents the ring of Laurent polynomials
  $\F[z_1,\ldots,z_m,z_1^{-1},\ldots,z_m^{-1}]$. In doing so a
  convolutional code then corresponds to an $\F$-linear subspace
  of some $\tilde{\Scal}^n$, where $\tilde{\Scal}=\{f:\Z^m
  \longrightarrow \F \mid f \text{ has finite support}\}$.
\end{rem}

Since $\D$ is a Noetherian ring, each code $\C\subseteq\D^n$ is
finitely generated.  In other words, there exists some $l\in\N$
and a matrix $G\in\D^{n\times l}$ such that $\C=\im_{\D}G$.  We
call such a matrix $G$ a generator matrix of $\C$.  Note that we
don't use the row vector notation as common in coding theory.  It
would force us to use the same notation also for the dual system
theoretic version, which is very unusual.  The notation
$\im_{\D}G$ means of course the set of all $Gp$ with $p\in\D^l$.
This notation instead of only $\im G$ will be necessary later
when interpreting $G$ as a different type of operator.
Analogously, we might also use the notation
$\ker_{\D}G=\{p\in\D^l\mid Gp=0\}$.

As a finitely generated $\D$-module each code has a well-defined
rank, say $\rank\C=k$.  It can simply be calculated as $\rank G$, where
one may use any generator matrix $G$ of $\C$, considered as a
matrix over the quotient field $\F(z_1,\ldots, z_m)$.  The rate
of $\C$ is defined to be the quotient $\frac{k}{n}$.

The code $\C$ is called free if $\C$ is a free $\D$-module, that
is, if $\C$ has a $\D$-basis.  This is the case if and only if
$\C$ has a generator matrix $G\in\D^{n\times k}$ with $\rank
G=k=\rank \C$.  Such a generator matrix is called an encoder.  If
$\C$ has an encoder, say $G=[G_1,\ldots,G_k]\in\D^{n\times k}$,
then each codeword can be written in a unique way as a
$\D$-linear combination of $G_1,\ldots,G_k$.  This is certainly a
very desirable property for a code.  It is a well-known fact that
each 1-dimensional code, that is each $\F[z_1]$-module is free.
However, for higher dimensions, i.~e. for $m>1$, this is not true
anymore.

\begin{exam}        \label{L0192}
  Let $\D=\F[z_1,z_2]$ and
\[
G(z_1,z_2)=
  \begin{bmatrix} 
    {z_1}^2& z_1+z_1 z_2\\
    z_1&1+z_2 \\
    z_1 z_2& z_2+{z_2}^2
  \end{bmatrix} 
\]
It can easily be shown that $\im_{\D}(G)$ is rate $\frac{1}{3}$
but not free.  That is, $\im_{\D}(G)$ has no $3 \times 1$
encoder.  The code $\im_{\D}\begin{bmatrix} z_1 & 1 & z_2
\end{bmatrix}^{\sf T}$ is free of rate $\frac{1}{3}$ and properly
contains $\im_{\D}(G)$.
\end{exam}

It is easy to see that encoder matrices for a given code are
unique up to unimodular right multiplication, i.~e. for
$G_i\in\D^{n\times k}$ with $\rank G_i=k$ it is
\begin{equation}\label{e-unique}
    \im_{\D}G_1=\im_{\D}G_2\Longleftrightarrow 
   G_2=G_1 U \text{ for some } U\in Gl_{k}(\D).
\end{equation}

An important measure for the `goodness' of a (convolutional) code
is its distance. In the remainder of this section we will
introduce this parameter.

\begin{de}     \label{L0224}
  Let $a \in \F^n$.  The {\it weight} of $a$ is given by the
  number of nonzero entries of $a$.  It is denoted by wt$(a)$.
  
  
  For $w=\sum_{\alpha \in \N^m} b_{\alpha} z^{\alpha}\in\D^n$ with
  $b_{\alpha} \in \F^n$ the {\it weight} of $w$ is defined as 
  $$
  \mbox{wt}(w)=\sum_{\alpha \in \N^m} \mbox{wt}(b_{\alpha})\,.
  $$
\end{de}

Hence the weight of a vector in $\D^n$ measures the distance to the all 
zero vector by counting all non-zero terms in 
the vector. The weight has the characteristic of a discrete
norm. In particular the weight induces a metric on $\D^n$:

\begin{de}     \label{L0232}
  Given two elements $w,\tilde{w} \in \D^n$ the {\it (Hamming)
    distance} between $w$ and $\tilde{w}$ is given by
  dist$(w,\tilde{w})=\mbox{wt}(w-\tilde{w})$.  Given any $m$-$D$
  code $\C$ of length $n$, the
  {\it distance of $\C$} is defined as
  $$
  \mbox{dist}(\C)={\rm min}\{{\rm dist}(w,\tilde{w}):w,
  \tilde{w} \in \C, w \neq \tilde{w}\}.
  $$
\end{de}

\begin{rem}             \label{L0240}
\begin{itemize}
\item[i)] The Hamming distance defines a metric on $\D^n$ called
  the Hamming metric.
\item[ii)] For a code $\C$ we have that
  dist$({\C})={\rm min}\{{\rm wt}(w):w \in {\C}, w \neq 0\}$.
  This is because dist$(w_1,w_2)={\rm wt}(w_1-w_2)$ and $w_1-w_2
  \in {\C}$ whenever $w_1,w_2 \in {\C}$.
\end{itemize}
\end{rem}

The following result is standard in coding theory. It follows
immediately from the definition of Hamming distance and from the
properties of a metric space.

\begin{lem}              \label{L0248} 
  Let $\C$ be a convolutional code with $d={\rm dist} ({\C}).$
  Let $t=\lfloor {\frac{d-1}{2}} \rfloor$ where $\lfloor x
  \rfloor$ denotes the greatest integer that is less than or
  equal to $x$.  Let $y \in \D^n$.  If $w \in {\C}$ is a codeword
  such that {\rm dist}$(w,y) \leq t$, then $w$ is the unique
  codeword nearest (with respect to the Hamming metric) to $y$.
\end{lem}

We say that $\C$ can correct up to $t$ errors. 
In practice it is not a simple task to compute the transmitted vector 
$w\in\D^n$ from the received vector $y\in\D^n$.
One way to do this is by syndrome decoding.
In order to explain this we first state

\begin{prop}           \label{L0136}
  Suppose $\C$ is a free convolutional code of rate $\frac{k}{n}$
  with encoder $G~\in~\D^{n \times k}$.  Then
  $$
  0 \longrightarrow \D^k \stackrel{G \cdot}{\longrightarrow}
  \D^n \stackrel{\pi}{\longrightarrow} \D^n / {\mathcal C}
  \longrightarrow 0
  $$
  is a short exact sequence.
\end{prop}

The proof is left to the reader. 

In the above short exact sequence $\pi$ is often called a syndrome former. 
Syndrome decoding works as follows.
If $y\in\D^n$ is received, one seeks the vector 
$e\in\pi^{-1}\big(\pi(y)\big)=y+\C$ of smallest possible weight.
The vector $y$ is then decoded as $w:=y-e$.

\section{Duality between Codes and Behaviors}
\setcounter{equation}{0}

There have been several instances in the recent literature about
coding theory, in which certain types of duality between
convolutional codes and behaviors in the system theoretical sense
of~\cite{wi91} have been mentioned or used, see,
e.~g.~\cite{fo95p3,ro96a1,va94,we98t}.

In this section we are going to make this duality precise by
introducing the appropriate bilinear form.  Exploiting the very
comprehensive and powerful paper of Oberst~\cite{ob90}, quite a
lot of results about this duality are available even in the
multidimensional case.  However, as it seems to us, most
interesting is the duality in the 1-dimensional case, where
various minimal first-order representations exist and have been
studied systematically and exhaustively by Kuijper~\cite{ku94}.
They can be translated into corresponding descriptions for codes.
This will be studied in Section~4.

We will introduce the notations and the setting along the
lines of \cite{ob90}.  Only those results needed for our purposes
will be cited afterwards.

Throughout this section the field $\F$ need not be finite; the
results hold for any field.  First we have to define the
underlying setting for the behaviors.  Let
\[
\A:=\F[\![z]\!]=\Big\{\sum_{\alpha\in\N^m}
f_{\alpha}z^{\alpha}\,\Big|\, f_{\alpha}\in\F\Big\}
\]
be the set of power series in the $m$ variables $z_1,\ldots,z_m$
over $\F$.  On $\A$ we consider the backward shifts along the
$i$th axis followed by truncation; that is, for each
$i=1,\ldots,m$ define
\begin{equation}\label{e-Lidef}
\begin{split}
  L_i\,:\qquad\ \A \quad &\longrightarrow \quad \A \\[1ex]
  \sum_{\alpha\in\N^m} f_{\alpha} z^{\alpha} &\longmapsto
  \sum_{\alpha\in\N^m} f_{\alpha+e_i}z^{\alpha}
\end{split}
\end{equation}
(see also \cite[p.~15]{ob90}).  This action can also be expressed
in the following ways
\begin{equation}                         \label{e-Li}
  L_i\Big(\sum_{\alpha\in\N^m} f_{\alpha}z^{\alpha}\Big)
   = z_i^{-1}\Big(\sum_{\2stack{\alpha
  \in\N^m}{\!\!\!\!\!\!\alpha_i\not=0}}
                          f_{\alpha}z^{\alpha}\Big)
   =\Pi_+\Big(z_i^{-1}\sum_{\alpha\in\N^m} f_{\alpha}z^{\alpha}\Big)
\end{equation}
where $\Pi_+$ denotes the projection which cuts off the terms
with negative exponents.  Clearly, the operators $L_i$ are
$\F$-linear.  Moreover, $\A$ gets the structure of a $\D$-module
via the scalar multiplication
\[
p(z_1,\ldots,z_m)\cdot f := p(L_1,\ldots, L_m)(f) \in\A \text{
  for } p\in\D,\, f\in\A.
\]

\begin{exam}\label{E-Li}
  Let $m=2$ and $p=1+z_1^2+z_1z_2\in\D=\F_5[z_1,z_2]$.  Then
  $p\cdot (1+3z_1z_2^3+2z_1^2+4z_2)=
  3+3z_1z_2^3+2z_1^2+4z_2+3z_2^2$.
\end{exam}

The example shows that the notation $p\cdot f$ has to be read
with care.  It is not the usual convolutional product in $\A$.
Instead from\eqr{e-Li} one can derive the formula
\begin{equation}\label{e-pdotf}
  \sideset{}{'}\sum_{\alpha\in\N^m}p_{\alpha}z^{\alpha}\cdot
     \sum_{\beta\in\N^m}f_{\beta}z^{\beta}
  = \sum_{\beta\in\N^m}
         \Big(\sideset{}{'}\sum_{\alpha\in\N^m} 
   p_{\alpha}f_{\alpha+\beta}\Big)z^{\beta}.
\end{equation}
Since we never use ordinary convolution in $\A$, this should not
cause a confusion.

\begin{rem}\label{R-Dmodule}
\begin{itemize}
\item[(a)] Obviously, $\D$ is a $\D$-submodule of $\A$.  However,
  it is worth mentioning that the canonical injection
  $\iota:\,\D\rightarrow \A, \;p\mapsto p$ is not $\D$-linear.
  In fact, e.~g., $z_1=\iota(z_1)\not=z_1\cdot\iota(1)=z_1\cdot
  1=0$ in $\A$.  This is not really an issue as the inclusion
  $\D\subset\A$ is never considered in this setting.  While $\D$
  is the set of operators, either generator matrices for codes or
  shift operators, $\A$ serves as the space of trajectories for
  the behaviors.
\item[(b)] $\A$ is not finitely generated as $\D$-module, see
  \cite[p.~55]{ob90}.
\end{itemize}
\end{rem}

Each polynomial matrix $G\in\D^{k\times n}$ gives rise to a
linear partial difference operator $G(L_1,\ldots,L_m)$ which we
will denote for short by $G$, thus
\[
\begin{split}
  G:\;\A^n&\longrightarrow \A^k\\[1ex]
  a &\longmapsto G\cdot a:=G(L_1,\ldots,L_m) (a)
\end{split}
\]
These operators are going to be the objects dual to generator
matrices for codes.  The following notations will be useful in
the sequel.  For $G\in\D^{k\times n}$ define
\[
\im_{\A}G=\{G\cdot a\mid a\in\A^n\}\subseteq
\A^k\quad\text{and}\quad \ker_{\A}G=\{a\in\A^n\mid G\cdot
a=0\}\subseteq\A^n.
\]
\begin{de}\label{D-behavior}
  An $m$-dimensional behavior $\B$ in $\A^n$ is defined to be a
  $\D$-submodule $\B\subseteq\A^n$ of the form $\B=\ker_{\A} G
  \text{ for some } G\in\D^{k\times n}$ (not necessarily of full
  row rank).
\end{de}

This setting is identical to the study of m-D-discrete-time
systems in the behavioral context, see e.~g.~\cite{ro90t}.

We observe that, while each $\D$-submodule of $\D^n$ is a code, not every 
$\D$-submodule of $\A^n$ is a behavior.
Characterizations for an $\F$-subspace of $\A^n$ being a behavior are given 
in the 1-dimensional case in \cite[III.1]{wi91} and for the general case
in \cite[p.~61/62]{ob90}. 

Now the bilinear form to be used for the duality is obvious.  For
each $n\geq1$ a $\D$-bilinear non-degenerate form is given by
(cf.~\cite[p.~22]{ob90})
\begin{equation}\label{e-bilinear}
\begin{split} 
  \D^n\times \A^n &\longrightarrow\quad  \A   \\
  (p, a)\quad& \longmapsto \dual{p}{a}:=p^{\sf T} \cdot a =
  \sum_{i=1}^n p_i\cdot a_i =\sum_{i=1}^n
  p_i(L_1,\ldots,L_m)(a_i)
\end{split}
\end{equation}
where $p=(p_1,\ldots,p_n)^{\sf T},\,a=(a_1,\ldots,a_n)^{\sf T}$.
In the literature related to codes and behaviors also a certain
$\F$-bilinear form has been used, see \cite{ro96a1} and
\cite[p.~20]{we98t}.  We will clarify the relationship between
this one and\eqr{e-bilinear} at the end of this section.

Using the above bilinear form we define the duals in the obvious
way.

\begin{de}\label{D-dual} 
\begin{itemize}
\item[(a)] The dual of a subset $\B\subseteq \A^n$ is defined to
  be $\B^{\perp}:=\{p\in\D^n\mid \dual{p}{a}=0\text{ for all }
  a\in\B\}$.
\item[(b)] The dual of a subset $\C\subseteq\D^n$ is given by
  $\C^{\perp}:=\{a\in\A^n\mid \dual{p}{a}=0\text{ for all }
  p\in\C\}$.
\end{itemize}
\end{de}

Obviously, duals are $\D$-modules and one has $\B\subseteq
\B^{\perp\perp}$ as well as $\C\subseteq \C^{\perp\perp}$.

Now we are in the position to state the results given
in~\cite{ob90}.  Essentially, they amount to the fact that $\A$ is
a large injective cogenerator in the category of $\D$-modules.
Instead of going into an explanation of this statement, we will
simply extract from \cite{ob90} the following consequences of
this very strong result.  Statements (4), (5), and  (7) of the next theorem are
exactly the duality between codes and behaviors we were looking
for.

\begin{thm}\label{T-lic}
  Let $P\in\D^{l\times n},\,Q\in\D^{k\times l},\,R\in\D^{r\times
    n}$. Then
\begin{itemize}
\item[\mbox{\rm(1)}] If the sequence $\D^k \stackrel{Q^{\sf
      T}}{\longrightarrow} \D^l \stackrel{P^{\sf
      T}}{\longrightarrow} \D^n$ is exact, then so is the
  sequence $\A^n \stackrel{P}{\longrightarrow} \A^l
  \stackrel{Q}{\longrightarrow} \A^k$.
\item[\mbox{\rm(2)}] $\ker_{\A}P\subseteq\ker_{\A}R$ if and only
  if $R=XP$ for some $X\in\D^{r\times l}$.
\item[\mbox{\rm(3)}] If $\rank P=l$, then the operator
  $P\,:\,\A^n\rightarrow\A^l$ is surjective.
\item[\mbox{\rm(4)}] $\big(\im_{\D}Q^{\sf T}\big)^{\perp}=
  \ker_{\A}Q$.
\item[\mbox{\rm(5)}] $\big(\ker_{\A}Q\big)^{\perp}=
  \im_{\D}Q^{\sf T}$.
\item[\mbox{\rm(6)}] $(\im_{\A}Q)^{\perp}=\ker_{\D}Q^{\sf T}$.
\item[\mbox{\rm(7)}] $\C=\C^{\perp\perp}$ and
  $\B=\B^{\perp\perp}$ for each code $\C\in\D^n$ and each
  behavior $\B\in\A^n$.
\end{itemize}
\end{thm}

(3) means in other words, for each $P$ with full row rank and for
each $g\in\A^l$ the associated linear partial difference equation
$P\cdot f=g$ has a solution in $\A^n$.  This is a well-known fact
in the 1-dimensional case, that is, $\D=\F[z_1]$.  Even more, one
can also prescribe initial conditions up to a certain order.
In the m-dimensional case this is more involved. 
Statement~(4) shows especially that the dual of a code is not
only a $\D$-module but even a behavior.

As for the proof, all the above results go back to
\cite[p.~33]{ob90}, which is just the large injective cogenerator
property.  However, we will give some more detailed references
and arguments from the paper to show how things are related with 
each other, although this might be a bit different from the order they have 
been proven.

(1) is exactly the injectivity of the module $\A$ which is
defined at \cite[p.~24]{ob90}.  (2) is at \cite[p.~36]{ob90}.
(3) is a consequence of (1).  (4) and (5) are at
\cite[p.~30/31]{ob90}, but they can also be derived directly from
the above as follows.  (4) and also (6) follow immediately from
\begin{equation}\label{e-dual}
   \dual{Q^{\sf T} p}{a}=\dual{p}{Q\cdot  a}\text{ for each }
   p\in\D^k\text{ and }a\in\A^l
\end{equation}
together with the non-degeneracy of the bilinear
form\eqr{e-bilinear}.  (5) can be shown with the help of (2) via
\[
\begin{split}
  p\in(\ker_{\A}Q)^{\perp} \Longleftrightarrow
  \ker_{\A}Q\subseteq\ker_{\A}p^{\sf T}%
  &\Longleftrightarrow p^{\sf T}=v^{\sf T} Q
  \text{ for some } v\in\D^k \\
  &\Longleftrightarrow p\in\im_{\D}Q^{\sf T}.
\end{split}
\]
(7) is a consequence from (4) and~(5).

\begin{rem} 
  (Compare with Remark~\ref{rem1}).  If convolutional codes are
  defined as submodules of $\tilde{\D}^n$, where $\tilde{\D}$
  represents the ring of Laurent polynomials $\F[z,z^{-1}]$ then
  this results in a duality between codes and linear behaviors
  defined on $\tilde{\A}^n$, where
  $\tilde{\A}:=\F[\![z,z^{-1}]\!]$ is the ring of formal power
  series in the variables
  $z_1,\ldots,z_m, z_1^{-1},\ldots,z_m^{-1}$.
\end{rem}

Next we want to concentrate on two specific descriptions of
behaviors.  They will be of significance for 1-dimensional
first-order-representations in the next section.  In fact, the
following two types of representations, applicable to both, codes
and behaviors, are dual to each other as will be proven next.
They specialize to the so-called $(P,Q,R)$- and
$(K,L,M)$-representations in the 1-dimensional case.
\begin{thm}\label{T-dual1}
  Let $R\in\D^{n\times l},\,N\in\D^{k\times l}$, and
  $M\in\D^{k\times n}$.  Then the following are true.
\begin{itemize}
\item[\mbox{\rm(a)}] The module
  $R\cdot(\ker_{\A}N):=\{R\cdot\zeta\mid
  \zeta\in\A^l,\,N\cdot\zeta=0\}\subseteq\A^n$ is a behavior and
  its dual is given by $\big(R\cdot(\ker_{\A}N)\big)^{\perp}
  =\{p\in\D^n\mid R^{\sf T} p\in\im_{\D}N^{\sf T}\}$.
\item[\mbox{\rm(b)}] The module $\{a\in\A^n\mid M\cdot
  a\in\im_{\A}N\}\subseteq\A^n$ is a behavior.  Its dual is
  $\{a\in\A^n\mid M\cdot a\in\im_{\A}N\}^{\perp}=
  M^{\sf T}(\ker_{\D}N^{\sf T})$.
\item[\mbox{\rm(c)}]
  $\big(R(\ker_{\D}N)\big)^{\perp}=\{a\in\A^n\mid R^{\sf T}\cdot
  a\in\im_{\A}N^{\sf T}\}$.
\item[\mbox{\rm(d)}] $\{p\in\D^n\mid M
  p\in\im_{\D}N\}^{\perp}=M^{\sf T}\cdot(\ker_{\A}N^{\sf T})$.
\end{itemize}
\end{thm}
\begin{proof}
  (a) The first part is proven in \cite[p.~26]{ob90}.  As for the
  second part, note the following equivalences, which hold for
  each $p\in\D^n$ using equation\eqr{e-dual}
\[
\begin{split}
\dual{p}{R\cdot a}=0\,\forall\,a\in\ker_{\A}N& \Longleftrightarrow
\dual{R^{\sf T} p}{a}=0\,\forall\,a\in\ker_{\A}N\\
&\Longleftrightarrow R^{\sf T}
p\in(\ker_{\A}N)^{\perp}=\im_{\D}N^{\sf T}.
\end{split}
\]
(b) Using $N=0$ in (a) we obtain especially that a $\D$-submodule
of the form $\im_{\A}R$ is a behavior.  Thus, write
$\im_{\A}N=\ker_{\A}Q$ with some appropriate $Q\in\D^{q\times
  k}$.  Then $\{a\in\A^n\mid M\cdot a\in\im_{\A}N\}=\ker_{\A}QM$
is a behavior (see also \cite[p.~27]{ob90}) and moreover
\[
\begin{split}
  \{a\in\A^n\mid M\cdot a\in\im_{\A}N\}^{\perp}
  &=(\ker_{\A}QM)^{\perp}=
  \im_{\D}(QM)^{\sf T}=M^{\sf T}(\im_{\D}Q^{\sf T})\\[1ex]
  &=M^{\sf T}\cdot\big((\ker_{\A}Q)^{\perp}\big)=M^{\sf
    T}\big((\im_{\A}N)^{\perp}\big) = 
    M^{\sf T}(\ker_{\D}N^{\sf T}).
\end{split}
\]
(c) and (d) follow now from (a) and (b) with
Thm.~\ref{T-lic}~(7).
\end{proof}

In the following we want to briefly discuss parity check matrices
for multidimensional codes.

\begin{de}\label{D-parity}
  Let $\C\subseteq\D^n$ be a code.  A matrix $H\in\D^{l\times n}$
  is called a parity check matrix of $\C$ if $\C=\ker_{\D}H$.
\end{de}

Not each code has a parity check matrix; e.~g. for $\D=\F[z_1]$
the code $\im_{\D}\Smallmat{z_1}{z_1}$ has no parity check
matrix, since each matrix $H\in\F[z_1]^{l\times2}$ having
$(z_1,z_1)^{\sf T}$ in its kernel, would also have $(1,1)^{\sf
  T}\in\ker_{\D}H$.

The following result about the existence of parity check matrices
can be found in \cite[3.3.8]{we98t}.

\begin{thm}\label{T-parity}
  Let $\C=\im_{\D}G$ with $G\in\D^{n\times k}$ be a free code,
  thus $\rank G=k$.  Then $\C$ has a parity check matrix if and
  only if $G$ is minor-prime, that is, if the greatest common
  divisor of all full-size minors of $G$ is a unit in $\D$.  If a
  parity check matrix exists, then one also has a parity check
  matrix $H\in\D^{(n-k)\times n}$ with $\rank H=n-k$.
\end{thm}

This result can be dualized by use of Thm.~\ref{T-lic}.

\begin{thm}\label{T-dual2}
  Let $\C\subseteq\D^n$ be a free code.  Then $\C$ has a parity check
  matrix if and only if the behavior $\C^{\perp}\subseteq\A^n$
  has an image-representation, i.~e.
\[
\C=\ker_{\D}H \text{ for some }H\in\D^{l\times
  n}\Longleftrightarrow \C^{\perp}=\im_{\A}H^{\sf T} \text{ for
  some }H\in\D^{l\times n}.
\]
Hence a behavior $\ker_{\A}G\subseteq\A^n$ has an image-representation
if and only if $G$ is minor-prime.
\end{thm}
\begin{proof}
 follows from Thm.~\ref{T-lic}~(1), (6),  and~(7).
\end{proof}

Recall that for 1-dimensional behaviors the existence of
image-representations is equivalent to controllability,
see~\cite{wi91}.  For $m>1$, at least one direction is true,
namely, behaviors with image-representations are always
controllable, see~\cite[Thm.~4.2]{wo97b}.  Equivalence can be
established for $m=2$ or for $m\geq2$ if certain directions of
the time-space axes are two-sided, see~\cite{ro90t} and
\cite[Thm.~6]{wo98}.

At the end of this section we want to discuss the relationship of
the above bilinear form with an $\F$-bilinear form which has
been used as well in the literature within this context.  Let
\[
\begin{split}
  \D^n\times \A^n \qquad\quad &\longrightarrow \qquad\F\\[1ex]
  \Big(\sideset{}{'}\sum_{\alpha\in\N^m} p_{\alpha}z^{\alpha},
  \sum_{\alpha\in\N^m}f_{\alpha}z^{\alpha}\Big) & \longmapsto
  \ddual{p}{f}:=\sideset{}{'}\sum_{\alpha\in\N^m} 
   p_{\alpha}^{\sf T} f_{\alpha}
\end{split}
\]
where $p_{\alpha}^{\sf T} f_{\alpha}\in\F$ denotes the usual
scalar product in $\F^n$.  Observe that the sum on the right hand
side is indeed finite.

\begin{exam}\label{E-ddual}
  Let $\D=\F_2[z_1]$ and $n=2$. For $p=(1,0)^{\sf T}\in\D^2$ and
  $f=(z_1,1)^{\sf T}\in\A^2$ we obtain
  $\ddual{p}{f}=(1,0)\Smallpmat{0}{1}=0$, whereas the previously
  used $\D$-bilinear form yields
  $\dual{p}{f}=(1,0)\Smallpmat{z_1}{1}=z_1$.  Hence $p$ and $f$
  are orthogonal with respect to $\ddual{\ }{\ }$ but not with
  respect to $\dual{\ }{\ }$.
\end{exam}

However, there is a close relationship between these two forms as
we will derive next.  In order to do so, we use the notation
$L^{\alpha}:=L_1^{\alpha_1}\circ\ldots\circ L_m^{\alpha_m}$ for
$\alpha=(\alpha_1,\ldots,\alpha_m)\in\N^m$ and the shifts 
$L_i$ defined in\eqr{e-Lidef}.  Let
$p=\sum'_{\alpha\in\N^m} p_{\alpha}z^{\alpha}\in\D^n$ and
$f=\sum_{\alpha\in\N^m}f_{\alpha}z^{\alpha}\in\A^n$.  Firstly,
using the very definition\eqr{e-bilinear} and
equation\eqr{e-pdotf} one obtains
\[
\begin{split}
  \dual{p}{f}=0&\Longleftrightarrow
  \sideset{}{'}\sum_{\alpha\in\N^m} p_{\alpha}^{\sf
    T}L^{\alpha}(f)=0 \Longleftrightarrow
  \sideset{}{'}\sum_{\alpha\in\N^m} p_{\alpha}^{\sf T}
  \sum_{\beta\in\N^m}f_{\beta+\alpha}z^{\beta}=0\\[1ex]
  &\Longleftrightarrow \sideset{}{'}\sum_{\alpha\in\N^m}
  p_{\alpha}^{\sf T}f_{\beta+\alpha}=0 \;\forall\,\beta\in\N^m
  \Longleftrightarrow \ddual{p}{z^{\beta}\cdot
    f}=0\;\forall\,\beta\in\N^m.
\end{split}
\]
Secondly, it is
\[
\ddual{z^{\beta}p}{f}=
\ddual{\sideset{}{'}\sum_{\alpha\in\N^m}p_{\alpha}
  z^{\alpha+\beta}}{\sum_{\alpha\in\N^m}f_{\alpha}z^{\alpha}}
=\sideset{}{'}\sum_{\alpha\in\N^m}p_{\alpha}^{\sf  T}
  f_{\alpha+\beta} =\ddual{p}{z^{\beta}\cdot f}
\]
by virtue of\eqr{e-pdotf}.

These two observations lead to the fact that both bilinear forms
yield the same duals for $\D$-submodules of $\A^n$ or $\D^n$.
Indeed, if $\B\subseteq\A^n$ is a $\D$-submodule, then
\[
\{p\in\D^n\mid \ddual{p}{f}=0\;\forall\,f\in\B\}= \{p\in\D^n\mid
\ddual{p}{z^{\beta}\cdot
  f}=0\;\forall\,f\in\B\;\forall\,\beta\in\N^m\}
=\B^{\perp}.
\]
Similarly, for a $\D$-submodule $\C\subseteq\D^n$ one obtains
\[
\begin{split}
  \{f\in\A^n \mid \ddual{p}{f}=0\;\forall\,p\in\C\}&= \{f\in\A^n
  \mid \ddual{z^{\beta}p}{f}=0\;\forall\,
  p\in\C\;\forall\,\beta\in\N^m\}\\[1ex]
  &=\{f\in\A^n \mid \ddual{p}{z^{\beta}\cdot
    f}=0\;\forall\,p\in\C\;\forall\,\beta\in\N^m\}
  =\C^{\perp}.
\end{split}
\]

\section{First-Order Representations for 1-Dimensional Codes}
\setcounter{equation}{0}

In this last section we restrict to the 1-dimensional case, thus
$\D=\F[z]$ denotes the polynomial ring in one variable over $\F$
and each submodule $\C\in\D^n$ is a convolutional code in the
sense of, e.~g., \cite{pi88}.  Using the duality results from the
last section and certain well-studied first-order representations
for behaviors, we can derive analogous descriptions for codes
along with minimality and uniqueness results.

The main source for this section is the book \cite{ku94} about
behaviors.  Although \cite{ku94} deals with the field $\R$, it
can be checked that the results hold true for any field.

We need to introduce the following parameter, called degree, for
1-dimensional codes.  It is the analogue to the McMillan degree
or order of a system.  Let $\C=\im_{\D}G$ with $G\in\D^{n\times
  k}$ and $\rank G=k$, a non-restrictive assumption.  The degree
$\delta(\C)$ is defined to be the maximum degree of all $k\times
k$-minors of $G$. The degree is sometimes also called the
complexity of the code $\C$ (see~\cite[2.7]{pi88}) and it
corresponds to the McMillan degree of the associated behavior
under the duality studied in the last section, see Thm~\ref{T-lic}~(4) and 
\cite[p.~276]{wi91}.
Equation\eqr{e-unique} shows that the degree does not depend on
the choice of the encoder $G$.  A code of degree $\delta(\C)=0$
is in essence a block code.

\begin{thm}\label{T-PQR}
  Let $\C=\im_{\D}G$ with $G\in\D^{n\times k}$ be a rate
  $\frac{k}{n}$ code of degree $\delta(\C)=\delta>0$.
\begin{itemize}
\item[\mbox{\rm (a)}] There exist matrices
  $(P,Q,R)\in\F^{\delta\times(\delta+k)}\times
  \F^{\delta\times(\delta+k)}\times\F^{n\times(\delta+k)}$
  such that
\[
\C=R\big(\ker_{\D}(zP+Q)\big).
\]
Moreover,
\\
\mbox{\rm (i)\;\;} $\rank P=\delta$,
\\
\mbox{\rm (ii)\;} $\rank\Smallmat{P}{R}=\delta+k$,
\\
\mbox{\rm (iii)} $zP+Q\in\D^{\delta\times(\delta+k)}$ is
left-prime.
\item[\mbox{\rm (b)}] If $\C=R\big(\ker_{\D}(zP+Q)\big)=\tilde{R}
  \big(\ker_{\D}(z\tilde{P}+\tilde{Q})\big)$ with matrix triples
  $(P,Q,R)$ and $(\tilde{P},\tilde{Q},\tilde{R})$ being of the
  sizes as in (a), then
\[
(\tilde{P},\tilde{Q},\tilde{R})=(T^{-1}PS,T^{-1}QS,RS) \text{ for
  some }T\in Gl_{\delta}(\F)\text{ and }S\in Gl_{\delta+k}(\F).
\]
\end{itemize}
\end{thm}
\begin{proof}
  (a) By Thm.~\ref{T-lic}~(4) we have 
   $\C^{\perp}=\ker_{\A}G^{\sf T}$.  
  Without loss of generality we may assume that $G$ is
  column-reduced, that is, $\delta$ is the sum of the column
  degrees of $G$.  From \cite[5.17]{ku94} we obtain matrices
  $(K,L,M)\in
  \F^{(\delta+k)\times\delta}\times\F^{(\delta+k)\times\delta}
  \times\F^{(\delta+k)\times n}$ such that
  $\C^{\perp}=\{a\in\A^n\mid M\cdot a \in\im_{\A}(zK+L)\}$.
  Indeed, the parameter ord$\,(\Sigma)$ in \cite[p.~128]{ku94} is
  equal to the degree, cf.  \cite[3.11 and 2.22]{ku94}.
  Setting $(P,Q,R)=(K^{\sf T},L^{\sf T},M^{\sf T})$ and using
  Thm.~\ref{T-dual1}~(b) and Thm.~\ref{T-lic}~(7) we obtain the
  desired representation.  Furthermore, \cite[5.17]{ku94} shows
  that the triple $(K,L,M)$ is minimal with respect to row and
  column size of the matrix $K$ (or $L$).  Hence, use of
  \cite[4.32]{ku94} leads to (i) -- (iii).
  
  (b) follows from \cite[4.40]{ku94} and Thm.~\ref{T-dual1}~(c).
\end{proof}

In fact, the proof shows more.  The above given sizes of the
matrices $(P,Q,R)$ are minimal among all representations of this
type.  The minimality is equivalent to the properties (i) --
(iii).  An alternative direct proof, without using duality, is
given in the paper \cite{sm98p3}.

In exactly the same way we can derive so-called
$(K,L,M)$-representations for codes. For this we use \cite[5.10
and 4.3]{ku94} and dualize these representations using
Thm.~\ref{T-dual1}~(a) and~(d). This results in~\cite[Theorem~3.1
and Theorem~3.4]{ro96a1}:

\begin{thm}\label{T-KLM}
  Let $\C=\im_{\D}G$ with $G\in\D^{n\times k}$ be a rate
  $\frac{k}{n}$ code of degree $\delta(\C)=\delta>0$.
\begin{itemize}
\item[\mbox{\rm (a)}] There exist matrices
  $(K,L,M)\in\F^{(\delta+n-k)\times\delta}
  \times\F^{(\delta+n-k)\times\delta}
  \times\F^{(\delta+n-k)\times n}$ so that
\[
\C=\{p\in\D^n\mid Mp\in\im_{\D}(zK+L)\}.
\]
Moreover,
\\
\mbox{\rm (i)\;\;} $\rank K=\delta$,
\\
\mbox{\rm (ii)\;} $\rank[K,M]=\delta+n-k$,
\\
\mbox{\rm (iii)} $[zK+L\mid M]$ is left-prime over the polynomial
ring $\D$.
\item[\mbox{\rm (b)}] If $\C=\{p\in\D^n\mid Mp\in\im_{\D}(zK+L)\}
  =\{p\in\D^n\mid \tilde{M}p\in\im_{\D}(z\tilde{K}+\tilde{L})\}$ with
  matrix triples $(K,L,M)$ and $(\tilde{K},\tilde{L},\tilde{M})$ being of the 
  sizes as in (a),  then
  \[
   (\tilde{K},\tilde{L},\tilde{M})=(T^{-1}KS,T^{-1}LS,T^{-1}M)
   \text{ for some }T\in Gl_{\delta+n-k}(\F),\ S\in Gl_{\delta}(\F).
  \]
\end{itemize}
\end{thm}

Generalized first order representations as described in the above two 
theorems are very useful in the design of convolutional codes
with large distance and which can be encoded in an
efficient manner. We refer the interested reader
to~\cite{ro96a1,ro97r1}. 

\section*{Conclusion}

The paper did show that multidimensional convolutional codes are
powerful encoding devices for the transmission of data over a
noisy channel. Since these codes are dual objects to
multidimensional systems the algebraic theory of linear systems
can be fruitfully applied.

Diederich Hinrichsen, to whom this paper is dedicated,
contributed over the years significantly to algebraic systems
theory. As it happens often in research a contribution in one
area bears unexpected fruits in another research field. We
believe that the recent cross fertilization between coding theory
and systems theory is such an instance. 


\end{document}